\documentclass[11pt]{amsart}
\usepackage{amsmath,amssymb,amsfonts,url,mathptmx}
\numberwithin{equation}{section}

\newtheorem{thm}{Theorem}[section]

\newtheorem{rem}{Remark}[section]

\begin{document}
\title[Vanishing estimates for Toda system]{Vanishing estimates for fully bubbling solutions of $SU(n+1)$ Toda Systems at a singular source} 

\author{Lei Zhang}
\address{Department of Mathematics\\
        University of Florida\\
        358 Little Hall P.O.Box 118105\\
        Gainesville FL 32611-8105}
\email{leizhang@ufl.edu}

\date{\today}

\begin{abstract} For Gauss curvature equation (or more general Toda systems) defined on two dimensional spaces, the vanishing rate of certain curvature functions on blowup points is a key
estimate for numerous applications. However, if these equations have singular sources, very few vanishing estimates can be found. In this article we consider a Toda system with singular sources
defined on a Riemann surface and we prove a very surprising vanishing estimates and a reflection phenomenon for certain functions involving the Gauss curvature.
\end{abstract}

\maketitle

\section{Introduction}

Let $(M,g)$ be a compact Riemann surface whose volume is assumed to be $1$ for convenience. Let $v=(v_1,...,v_n)$ be a solution of the following $SU(n+1)$ Toda system defined on $M$:
\begin{equation}\label{eq-0}
\Delta_g  v_i+\sum_{j=1}^n a_{ij}H_je^{v_j}-K(x)=4\pi \sum_m \gamma_{im}\delta_{q_m},\quad 1\le i\le n
\end{equation}
where $\Delta_g$ is the Laplace-Beltrami operator ($-\Delta_g\ge 0$), $K$ is the Gauss curvature, $H_1,..,H_n$ are positive smooth functions on $M$, $\delta_{q_m}$ is the Dirac mass at $q_m$, $\gamma_{im}\in \mathbb R$ is assumed to be greater than $-1$ for all $i$ and $m$, $A=(a_{ij})_{n\times n}$ is the following Cartan matrix:
$$A=\left(\begin{array}{ccccc}
2 & -1 & 0 & ... & 0 \\
-1 & 2 & -1 & ... & 0 \\
\vdots & \vdots & \vdots & ... & \vdots \\
0& ... & -1 & 2 & -1 \\
0 & ... & ... & -1 & 2
\end{array}
\right )
$$
The $SU(n+1)$ Toda systems are well known to be deeply rooted in algebraic geometry (see \cite{calabi,doliwa,jostlinwang,lin-wei-ye,lin-wei-zhang-jems,lin-wei-zhang-adv,lin-wei-zhao}) and have close connections with various fields such as integrable system, the non-abelian Chern-Simons model in Gauge theory, etc (see \cite{dunne1,dunne2,ganoulis,guest,leznov1,leznov2,mansfield,nolasco1,nolasco2,prajapat,yang1,yang2} and Painleve VI equation\cite{chai,CKL}, etc.  The readers may read the introduction of the
first two articles in this series \cite{lin-wei-zhang-jems,lin-wei-zhang-adv} and \cite{chai,CKL} for the references and explanations in more detail. The purpose of this article is to study the behavior of blowup solutions of the singular Toda system if the blowup point happens to be a singular source.

Equation (\ref{eq-0}) is usually written in the following form:
\begin{equation}\label{toda-1}
\Delta_g u_i+\sum_{i=1}^na_{ij}\rho_j(\frac{h_je^{u_j}}{\int_M h_je^{u_j}dV_g}-1)=\sum_{m=1}^N 4\pi\gamma_{im}(\delta_{q_m}-1),\quad i=1,..,n
\end{equation}
where $A=(a_{ij})_{n\times n}$ is the Cartan matrix, $h_1,...,h_n$ are positive smooth functions on $M$, $\rho_1,..,\rho_n$ are positive constants, $q_1$, ..., $q_N$ are distinct points on $M$ and $\gamma_{im}$ are all greater than $-1$.
In this article we suppose $q_1$ is a blowup point, (note that at $q_1$ there is a Dirac source) and we consider the behavior of blowup solutions in the neighborhood of $q_1$.  Since our main result is local in nature, we state our result for the following locally defined $SU(n+1)$ Toda system for convenience:

Let $u_k=(u_{k,1},...,u_{k,n})$ be a sequence of solutions of
\begin{equation}\label{eq-1}
\Delta u_{k,i}+\sum_{j=1}^n a_{ij}h_{k,j}(x)e^{u_{k,j}}=4\pi \gamma_i\delta_0,\quad \mbox{ in }\quad B_1,
\quad i=1,..,n,\quad n\ge 2
\end{equation}
where $B_1$ is the unit ball in $\mathbb R^2$ ( throughout the article we use $B(x_0,r)$ to denote the ball centered at $x_0$ with radius $r$. If $x_0$ is the origin we shall just use $B_r$), $A=(a_{ij})_{n\times n}$ is the Cartan matrix,
 $\delta_0$ is the Dirac source at the origin,
$h_{k,i}$ are positive and smooth functions in $B_1$: There exists $c>0$ such that
\begin{equation}\label{hik}
\frac 1c\le h_{k,i}(x)\le c,\quad \forall x\in B_1, \quad |D^2h_{k,i}(x)|\le c,\quad \forall x\in B_1.
\end{equation}
For $\gamma_i$ we assume
\begin{equation}\label{gamma-i2}
\gamma_i\ge 0 \quad \mbox{ for all } i=1,..,n.
\end{equation}
and we shall use $I_1$ to denote the collection of nonzero indexes and $I_2$ to denote the set of zero indexes:
$$I_1=\{i;\quad \gamma_i\neq 0\}, \quad I_2=\{i;\quad \gamma_i=0\} $$

In addition we make the following natural assumption, which is pretty much postulated in all the works related to the study of blowup solutions in two dimensional spaces:
\begin{equation}\label{assumption-m}
\left\{\begin{array}{ll}
\max_{k\subset\subset B_1\setminus \{0\}} u_{k,i} \le C(K),\quad i=1,...,n,\quad n\ge 2, \\
\\
\max_{x,y\in \partial B_1}|u_{k,i}(x)-u_{k,i}(y)|\le C.  \\
\\
\int_{B_1}h_{k,i}e^{u_{k,i}}\le C \mbox{ for some $C>0$ independent of $k$ }.
\end{array}
\right.
\end{equation}

Since $u_k=(u_{k,1},..,u_{k,n})$ has a logarithmic term corresponding to the singular source, it is convenient to consider the equation for the regular part of $u_k$. Let
$\tilde u_k=(\tilde u_{k,1},...,\tilde u_{k,n})$ be the regular part of $u_k$:
\begin{equation}\label{tuk}
\tilde u_{k,i}(x)=u_{k,i}(x)-2\gamma_i\log |x|,\quad i=1,..,n.
\end{equation}
Then we have
\begin{equation}\label{eq-2}
\Delta \tilde u_{k,i}+\sum_{j=1}^n a_{ij} h_{k,j}(x) |x|^{2\gamma_j}e^{\tilde u_{k,j}(x)}=0,\quad \mbox{in }\quad B_1.
\end{equation}
Since we study blowup solutions, the maximum of $\tilde u_k$ tends to infinity: Let
$$M_k=\max_i\max_{x\in B_1}\frac{|\tilde u_{k,i}(x)|}{1+\gamma_i}\to \infty,\quad \epsilon_k=e^{-\frac 12 M_k}, $$
and
\begin{equation}\label{vik-d}
\tilde v_{k,i}(y)=\tilde u_{k,i}(\epsilon_ky)+2(1+\gamma_i)\log \epsilon_k,\quad i=1,...,n.
\end{equation}
Then direct computation shows
\begin{equation}\label{vik}
\Delta \tilde v_{k,i}(y)+\sum_j a_{ij}h_{k,j}(\epsilon_k y)|y|^{2\gamma_j}e^{\tilde v_{k,j}(y)}=0,\quad |y|\le \epsilon_k^{-1}.
\end{equation}
It is well known that for systems, if the whole system is scaled according to the maximum of all components, it is possible to have some components tending to minus infinity over any fixed compact subsets, which means these components do not appear in the limiting system. Such a situation is called partial blowup phenomenon. If no component is lost in the limit system, such a blowup sequence is called fully bubbling. The main assumption in this article is that $\tilde v_k=(\tilde v_{k,1},...,\tilde v_{k,n})$ is fully bubbling :
\begin{equation}\label{maj-a}
\tilde v_k=(\tilde v_{k,1},...,\tilde v_{k,n}) \to \tilde v=(\tilde v_1,...,\tilde v_n) \mbox{ in } C^{1,\alpha}_{\mbox{loc}}(\mathbb R^2) \quad \alpha\in (0,1)
\end{equation}
where $\tilde v=(\tilde v_1,...,\tilde v_n)$ satisfies
\begin{equation}\label{fully-bu}
 \left\{\begin{array}{ll}
\Delta \tilde v_i+\sum_{j=1}^n a_{ij} |y|^{2\gamma_j}e^{\tilde v_j}=0,\quad \mbox{in}\quad \mathbb R^2,\quad i=1,..,n \\
\\
\int_{\mathbb R^2}|y|^{2\gamma_i}e^{\tilde v_i}dy<\infty, \quad i=1,..,n.
\end{array}
\right.
\end{equation}
Here for convenience we assume $\lim_{k\to \infty} h_i^k(0)=1$, but this assumption is not essential.

Let $\psi_{k,i}$ be a harmonic function that makes $u_{k,i}-\psi_{k,i}$ constant on $\partial B_1$:
\begin{equation}\label{psi-k}
\left\{\begin{array}{ll}
\Delta \psi_{k,i}=0,\quad \mbox{in}\quad B_1,\\
\psi_{k,i}(x)=u_{k,i}(x)-\frac 1{2\pi}\int_{\partial B_1}u_{k,i} \quad \mbox{ on } \partial B_1.
\end{array}
\right.
\end{equation}
Since $u_{k,i}$ has bounded oscillation on $\partial B_1$ we have
$$\psi_{k,i}(0)=0,\quad |D^m \psi_{k,i}|_{L^{\infty}(B_{1/2})}\le C(m),\quad \mbox{ for } m=0,1,2,.. $$

The main result of this article is
\begin{thm}\label{local-est} Suppose $I_2$ is not empty and $n\ge 2$, $\tilde v_k$ is a fully bubbling sequence described in (\ref{maj-a}), $h_k$ and $\gamma=(\gamma_1,...,\gamma_n)$ satisfy the conditions stated in
(\ref{hik}),(\ref{gamma-i2}) and (\ref{assumption-m}). $\psi_k$ is defined in (\ref{psi-k}). Then we have
$$|\nabla (\log h_{k,i}+ \psi_{k,i})(0)|=O(\epsilon_k),\quad \mbox{ if } n+1-i\in I_2. $$
\end{thm}

If the Toda system has only one equation, it is a prescribing Gauss curvature equation which has been extensively studied for decades. It is well known that the location of a blowup point has to be a critical point of a curvature function and this function has to vanish to $0$ along the local maximums of blowup solutions at a certain rate. Such a key estimate has a number of important applications such as constructing the blowup solutions, capturing the geometric information of the manifold, computing the Leray-Schauder degree of the set of solutions in terms of the topology of the manifold, etc. However, almost no vanishing result can be found if the blowup point happens to be a singular source. In other words, even for the singular equation, the study of the vanishing rate of curvature function at a singular source seems completely blank. To the best of our knowledge,
Theorem \ref{local-est} provides the first vanishing estimate for singular Toda systems. Moreover, the reflection of an index (from $i$ to $n+1-i$) is also a new feature that has never been observed before.

For the purpose of application Theorem \ref{local-est} can be written under a more general setting:

Let $\mu_i^k$ be a sequence of smooth functions tending to
$4\pi \gamma_i\delta_0$ in measure. We assume that $u^k=(u_1^k,..,u_n^k)$ satisfies (\ref{assumption-m}) and
\begin{equation}\label{gen-e1}
\Delta u_i^k+\sum_{j=1}^n a_{ij} h_j^k(x) e^{u_j^k}=\mu_i^k, \quad \mbox{ in }\quad B_1,
\end{equation}
where $h^k=(h_1^k,...,h_n^k)$ satisfies (\ref{hik}) as well.
Then we use $f_i^k$ defined in the following as a replacement of $2\gamma_i\log |x|$:
\begin{equation}\label{gen-2}
\left\{\begin{array}{ll}
\Delta f_i^k=\mu_i^k,\quad  \mbox{ in }\quad B_1, \\
f_i^k=0,\quad \mbox{on }\quad \partial B_1.
\end{array}
\right.
\end{equation}
Now the definition of $\tilde u^k=(\tilde u_1^k,...,\tilde u_n^k)$ becomes
$$\tilde u_i^k(x)=u_i^k(x)-f_i^k. $$
Clearly $\tilde u^k$ satisfies
\begin{equation}\label{gen-tuk}
\Delta \tilde u_i^k+\sum_j a_{ij} h_j^ke^{f_j^k}e^{\tilde u_j^k}=0,\quad \mbox{ in }\quad B_1,\quad i=1,..,n.
\end{equation}
Let
\begin{equation}\label{ea-1}
M_k=\max_i\max_{x\in B_1}\frac{\tilde u_i^k(x)}{1+\gamma_{i,k}}
\end{equation}
where $\gamma_{i,k}$ is a sequence of constants tending to $\gamma_i$. The specific requirements of $\gamma_{i,k}$ will be stated later.
Set
\begin{equation}\label{ea-2}
\epsilon_k=e^{-\frac 12 M_k}
\end{equation}
and
\begin{equation}\label{ea-3}
\tilde v_i^k(y)=\tilde u_i^k(\epsilon_k y)+2(1+\gamma_{i,k})\log \epsilon_k, \quad i=1,..,n.
\end{equation}
Direct computation shows
\begin{equation}\label{gen-tvk}
\Delta \tilde v_i^k+\sum_ja_{ij}h_j^k(\epsilon_k y)\frac{e^{f_j^k(\epsilon_ky)}}{\epsilon_k^{2\gamma_{i,k}}}e^{\tilde v_j^k(y)}=0,\quad |y|\le 1/\epsilon_k^{-1}
\end{equation}
It is easy to see from the definition of $\epsilon_k$ that
$$\tilde v_i^k\le 0 \mbox{ in  }B(0,\epsilon_k^{-1})\quad i=1,..,n, \mbox{ and  }
\max_i\max_y \tilde v_i^k(y)=0, $$
 so our main assumption in this more general setting is
\begin{align}\label{gen-a}
\mbox{ there exists }\gamma_{i,k}\to \gamma_i \mbox{ such that } \tilde v^k=(\tilde v_1^k,...,\tilde v_n^k) \mbox{ converges uniformly to }\\
(\ref{fully-bu}) \mbox{ over any fixed compact subset of } \mathbb R^2. \nonumber
\end{align}
Then the conclusion of Theorem \ref{local-est} also holds for $u^k$ in this more general setting:

\medskip

\begin{thm}\label{thm-general-v}: Let $u^k=(u_1^k,...,u_n^k)$ be a sequence of solutions to (\ref{gen-e1}) that satisfies (\ref{assumption-m}), $\mu^k=(\mu_1^k,...,\mu_n^k)$ be a sequence of smooth functions tending to $4\pi \gamma_i\delta_0$ in measure (with $\gamma=(\gamma_1,...,\gamma_n)$ satisfying (\ref{gamma-i2})). Let $f_i^k$, $M_k$, $\epsilon_k$ be defined by (\ref{gen-2}), (\ref{ea-1}) and (\ref{ea-2}), respectively. Then for $\tilde v^k=(\tilde v_1^k,...,\tilde v_n^k)$ in (\ref{ea-3}), if (\ref{gen-a}) holds,
$$|\nabla (\log h_i^k+\psi_{k,i})(0)|=O(\epsilon_k), \mbox{ if  } n+1-i\in I_2 $$
where $\psi_{k,i}$ is defined in (\ref{psi-k}).
\end{thm}

Even though Theorem \ref{local-est} and Theorem \ref{thm-general-v} are stated for locally defined Toda systems, they are very useful for equation (\ref{toda-1}).
One major goal for studying (\ref{toda-1}) is to identify the set of critical parameters when the blowup phenomenon occurs. This information is related to the topology of $M$ and a corresponding degree counting formula. In order to achieve this goal one major difficulty comes from the asymptotic behavior of blowup solutions near an isolated blowup point. In general the asymptotic behavior of blowup solutions is very complicated and becomes significantly more difficult as the number of equations increases. Theorem \ref{local-est} and Theorem \ref{thm-general-v} should be very useful for simplifying bubble interactions even for regular $SU(n+1)$ Toda systems. In fact, even for solutions to the regular $SU(n+1)$ Toda system, if different components tend to infinity at different speed, the components that tend to infinity fast look like Dirac mass for slower components. The phenomenon that Theorem \ref{local-est} and Theorem \ref{thm-general-v} revealed is quite new and no similar results have been observed or verified for similar systems before.

The sequence of harmonic function $\psi_k=(\psi_{k,1},...,\psi_{k,n})$ is usually easy to be identified in application. For example for blowup solutions of (\ref{toda-1}), if the blowup sequence $u_k$ converges in measure to a few Dirac masses and some rough estimates for $u_k$ are obtained outside the bubbling area, $\psi_k$ can be determined easily by the Green's function of $-\Delta_g$ on the manifold.

To indicate the application of Theorem \ref{local-est} we present the following simple example. Let $M$ be a flat torus and we consider
$$\left\{\begin{array}{ll}\Delta_g u_1+2e^{u_1}-e^{u_2}=4\pi \gamma_1 \delta_{p_1}+4\pi \gamma_2\delta_{p_2}\\
\Delta_g u_2-e^{u_1}+2e^{u_2}=4\pi \gamma_3 \delta_{p_3}.
\end{array}
\right.
$$
Here we assume $p_1$, $p_2$, $p_3$ are distinct points and $p_1$ is the only blowup point and the blowup sequence is fully bubbling at
$p_1$. Let $(u_{k,1},u_{k,2})$ be the fully bubbling sequence, then it is easy to use the classification theorem of Lin-Wei-Ye \cite{lin-wei-ye} to obtain $e^{u_{k,1}}\rightharpoonup 4\pi(2+\gamma_1)\delta_{p_1}$ and
$e^{u_{k,2}}\rightharpoonup 4\pi(2+\gamma_1)\delta_{p_1}$. Thus $\gamma_2=2$ and $\gamma_3=2+\gamma_1$.
The Green's representation formula for $u_1^k$ gives
$$u_{k,1}(x)=\bar u_{k,1}+\int_M G(x,\eta)(2 e^{u_{k,1}}-e^{u_{k,2}}-4\pi\gamma_1\delta_{p_1}-8\pi \delta_{p_2})dV_g(\eta) $$
where $\bar u_{k,1}$ is the average of $u_{k,1}$ and $G$ is the Green's function that satisfies
$$\Delta_xG(x,p)=-\delta_p+\frac 1{vol(M)},\quad \int_M G(x,p)dV_g(x)=0. $$
Using the concentration of $u_{k,1}$ and $u_{k,2}$ we see that in the neighborhood of $p_1$
$$ u_{k,1}(x) =\bar u_{k,1}+8\pi(G(x,p_1)-G(x,p_2))+o(1), \quad x\in B(p_1,2\delta)\setminus B(p_1,\delta) $$
for some small $\delta>0$. The vanishing estimate in Theorem \ref{local-est} gives
\begin{equation}\label{extra-ex}
\nabla_1\gamma(p_1,p_1)-\nabla_1G(p_1,p_2)=0.
\end{equation}
where $\gamma$ is the regular part of $G$ and $\nabla_1$ means the differentiation with respect to the first component. In other words
if $p_1$, $p_2$ do not satisfy (\ref{extra-ex}) it is not possible to have $p_1$ as the only blowup point with a fully bubbling sequence.

In \cite{lin-wei-zhang-adv} three sharp estimates are obtained for fully bubbling solutions of regular
$SU(n+1)$ Toda systems, the third of which is a $\partial_z^2$ condition: an estimate on the second derivatives of coefficient functions at blowup points. The readers may wonder why that estimate is not derived in this article. The reason is the estimate in Theorem \ref{local-est} is not for all the indexes (however the corresponding estimate in \cite{lin-wei-zhang-adv} holds for all indexes). This fact prevents us from getting more accurate asymptotic behavior of fully bubbling solutions and the $\partial_z^2$ conditions.

The proof of Theorem \ref{local-est} relies heavily on the important classification theroem of Lin-Wei-Ye \cite{lin-wei-ye} on global solutions of $SU(n+1)$ Toda system. In particular we find out that some leading terms in the asymptotic behavior of global solutions are crucial for blowup analysis. By differentiating on certain parameters related to these leading terms we obtain a few families of solutions to the linearized $SU(n+1)$ Toda systems, which play an important role in the proof of Theorem \ref{local-est}.

The organization of this article is as follows. In section two we study the asymptotic behavior of global solutions to singular $SU(n+1)$ Toda systems. It turns out that the components of the global solution that correspond to $I_2$ contribute some crucial leading terms in the expansion of the global solution. Then in section three we prove Theorem \ref{local-est}. The proof of Theorem \ref{thm-general-v} is very similar to that of Theorem \ref{local-est} and is therefore omitted.

\section{Properties of global solutions to $SU(n+1)$ Toda system with one singular point}
In this section we study the asymptotic behavior of global solutions to $SU(n+1)$ Toda system with one singularity. The main results of this section are (\ref{exp-u}) and (\ref{U-diff}).

\medskip

Let $U=(U_1,...,U_n)$ satisfy
\begin{equation}\label{toda-nk}
\left\{\begin{array}{ll}
\Delta U_i+\sum_{j=1}^n a_{ij}e^{U_j}=4\pi \gamma_i\delta_0,\quad i=1,...,n \quad \mbox{in}\quad \mathbb R^2,\\
\\
\int_{\mathbb R^2}e^{U_i}<\infty,
\end{array}
\right.
\end{equation}
where $A=(a_{ij})_{n\times n}$ is the Cartan matrix.
Recall that $I_1$ is the set of nonzero indexes and $I_2$ is the complement of $I_1$.

Let
$$U^i=\sum_ja^{ij}U_j,\quad \gamma^i=\sum_j a^{ij}\gamma_j. $$

We shall use the following properties of global solutions ( see \cite{lin-wei-ye} ):
$$e^{-U^1}=f, \quad e^{-U^k}=2^{k(k-1)}det_k(f) $$
where
\begin{align*}
&det_k(f)=det(f^{(p,q)}),\quad 0\le p,q\le k-1,\\
& f^{(p,q)}=\partial_{\bar z}^q\partial_z^pf, \\
&|z|^{2\gamma^1}e^{-U^1}=\lambda_0+\sum_{1\le i\le n}\lambda_i |P_i(z)|^2 \\
&P_i(z)=\sum_{j=0}^{i}c_{ij}z^{\mu_1+...+\mu_j},\,\, \mu_0=0,\,\, \mbox{and}\, \mu_i=1+\gamma_i, c_{ii}=1 \,\, \mbox{for} \,\, i=1,..,n,.
\end{align*}
Moreover
$$\lambda_0...\lambda_n=2^{-n(n+1)}\Pi_{1\le i\le j\le n}(\sum_{k=i}^j\mu_k)^{-2}, $$
and
$$f=\lambda_0|z|^{-2\gamma^1}+\sum_{i=1}^n\lambda_i P_i(z)z^{-\gamma^1}\bar P_i(z)\bar z^{-\gamma^1}. $$

Let  $q_0(z)=z^{-\gamma^1}$ and
$$q_i(z)=\sum_{j=0}^{i}c_{ij}z^{\mu_1+...+\mu_j-\gamma^1},\quad i=1,...,n $$
we have
$$f=\sum_{i=0}^n\lambda_iq_i(z)\bar q_i(z) $$
and
$$f^{(p,q)}=\sum_{i=0}^n\lambda_i \partial_p q_i\partial_q \bar q_i. $$
\begin{align*}
e^{-U^m}&=2^{m(m-1)}det \bigg ( (C_{\lambda},B_{\lambda})\left (\begin{array}{c} \bar C \\
\bar B
\end{array}
\right )  \bigg ) \\
&=2^{m(m-1)}det(C_{\lambda}\bar C+B_{\lambda}\bar B).
\end{align*}
where
$$B_{\lambda}=\left(\begin{array}{ccc}
\lambda_{n+1-m}q_{n+1-m} & ... & \lambda_n q_n \\
\lambda_{n+1-m}q_{n+1-m}^{(1)} & ... & \lambda_n q_n^{(1)}\\
\vdots &  \vdots & \vdots \\
\lambda_{n+1-m}q_{n+1-m}^{(m-1)} & ... & \lambda_n q_n^{(m-1)}
\end{array}
\right )
$$
$$C_{\lambda}=\left(\begin{array}{ccc}
\lambda_{0}q_{0} & ... & \lambda_{n-m} q_{n-m} \\
\lambda_{0}q_{0}^{(1)} & ... & \lambda_{n-m} q_{n-m}^{(1)}\\
\vdots &  \vdots & \vdots \\
\lambda_{0}q_{0}^{(m-1)} & ... & \lambda_{n-m} q_{n-m}^{(m-1)}
\end{array}
\right )
$$
$$
\bar B'=\left(\begin{array}{cccc}
\bar q_{n+1-m} & \bar q_{n+1-m}^{(1)} & ... & \bar q_{n+1-m}^{(m-1)} \\
\vdots & \vdots & \vdots & \vdots \\
\bar q_n & \bar q_n^{(1)} & ... & \bar q_n^{(m-1)}
\end{array}
\right )
$$
$$
\bar C'=\left(\begin{array}{cccc}
\bar q_{0} & \bar q_{0}^{(1)} & ... & \bar q_{0}^{(m-1)} \\
\vdots & \vdots & \vdots & \vdots \\
\bar q_{n-m} & \bar q_{n-m}^{(1)} & ... & \bar q_{n-m}^{(m-1)}
\end{array}
\right )
$$
where $q_i^{(j)}$ means the $j$-th $z$ derivative of $q_i$. $\bar q_i^{(j)}$ is the $j$th $\bar z$ derivative of $\bar q_i$.
The leading term comes from $B_{\lambda}\bar B'$. Our goal is to determine the first two terms in the expansion of $e^{-U^m}$. Here we note that $C_{\lambda}$ and $\bar C$ may not be square matrices, but $B_{\lambda}$ and $\bar B$ are
square matrices and $det(B_{\lambda}\bar B)$ will give at least the first two leading terms.
Let
$$s_i=\mu_1+...+\mu_i-\gamma^1,\quad i=1,..,n. $$
Using the fact that  $0<a^{ij}<1$ and the definition of $\gamma^1$ we see that $s_n>n$.  Next we see that
$$q_i=z^{s_i}+c_{i,i-1}z^{s_{i-1}}+c_{i,i-2}z^{s_{i-2}}+l.o.t $$
where $l.o.t$ stands for ``lower order terms". From the definition of $s_i$ we see that
$$s_{i-1}=s_i-1,\quad \mbox{ if }\quad \gamma_i=0; \quad s_i>s_{i-1}+1\quad \mbox{if}\quad \gamma_i>0. $$
Consequently if $i\in I_2$, $s_i=s_{i-1}+1$, and $s_i\ge s_{i-2}+2$. If $i\in I_1$ and $\gamma_i\not \in \mathbb N$ (the set of all natural numbers), $c_{i,i-1}=0$ 
(see the main theorem of \cite{lin-wei-ye}) and $s_{i-2}\le s_i-2$ still holds. If $i\in I_1$ and $\gamma_i$ is a positive integer, $s_i\ge s_{i-2}+2$ clearly holds. Therefore for $z$ large
\begin{align*}
&q_i=z^{s_i}(1+c_{i,i-1}/z+O(1/z^2), \quad \mbox{if }\quad i\in I_2, \\
&q_i=z^{s_i}(1+O(1/z^2)),\quad \mbox{if }\quad i\in I_1.
\end{align*}

In order to identify the two leading terms of $U^m$, we first identify the leading term in $B_{\lambda}$. By taking out $\lambda_{n+1-m},..,\lambda_n$ and ignoring all the l.o.t in each entry we have
\begin{align*}
&det B_{\lambda}=\lambda_{n+1-m}...\lambda_n \cdot\\
&det \left(\begin{array}{ccc}
z^{s_{n+1-m}} & ... & z^{s_n}\\
s_{n+1-m}z^{s_{n+1-m}-1} & ...  & s_n z^{s_n-1} \\
\vdots & \vdots & \vdots\\
\Pi_{j=0}^{m-2}(s_{n+1-m}-j)z^{s_{n+1-m}+1-m} & ... & \Pi_{j=0}^{m-2}(s_n-j)z^{s_n+1-m}
\end{array}
\right )+l.o.t.
\end{align*}
Note that the l.o.t are with respect to the leading term in the determinant.
For the major matrix we first take out $z^{s_{n+1-m}+1-m}$ from the first column, $z^{s_{n+2-m}+1-m}$ from the second column,..., $z^{s_n+1-m}$ from the $m-th$ column, then the power of entries in the first row becomes $z^{m-1}$, in the second row it is $z^{m-2}$,... in the $m-1$th row it is $z$ and in the last row there is no $z$. Thus by taking out $z^{m-1}$ from the first row, $z^{m-2}$ in the second row,.., $z$ from the $m-1$th row we see that the power of $z$ of the leading matrix is
$$S_m:=s_{n+1-m}+...+s_n-\frac{m(m-1)}2. $$
 Now $det B_{\lambda}$ becomes
$$det B_{\lambda}=\lambda_{n+1-m}...\lambda_nz^{S_m}det\left(\begin{array}{ccc}
1 & ... &1\\
s_{n+1-m}  & ...  & s_n   \\
\vdots & \vdots & \vdots\\
\Pi_{j=0}^{m-2}(s_{n+1-m}-j) & ... & \Pi_{j=0}^{m-2}(s_n-j)
\end{array}
\right )+l.o.t
$$
The evaluation of
$$D=det\left(\begin{array}{ccc}
1 & ... &1\\
s_{n+1-m}  & ...  & s_n   \\
\vdots & \vdots & \vdots\\
\Pi_{j=0}^{m-2}(s_{n+1-m}-j) & ... & \Pi_{j=0}^{m-2}(s_n-j)
\end{array}
\right ) $$
is elementary. First we subtract $(s_n+2-m)$ times row $m-1$ from row $m$. Then the entries of the last row are
$$\Pi_{j=0}^{m-3}(s_{n+1-m}-j)(s_{n+1-m}-s_n),...\Pi_{j=0}^{m-3}(s_{n-1}-j)(s_{n-1}-s_n),0. $$
Next we substract $s_n+3-m$ times row $m-2$ from row $m-1$. Then the $m-1$th row is the following after this operation:
$$\Pi_{j=0}^{m-4}(s_{n+1-m}-j)(s_{n+1-m}-s_n),...,0. $$
Eventually we substract $s_n$ times row $1$ from row $2$ and the second row becomes
$$s_{n+1-m}-s_n,...,s_{n-1}-s_n,0. $$
By expanding at the $(1,m)$ entry we see that the determinant is equal to
$$\Pi_{j=n+1-m}^{n-1}(s_n-s_j)
det\left(\begin{array}{ccc}
1 & ... &1\\
s_{n+1-m}  & ...  & s_{n-1}   \\
\vdots & \vdots & \vdots\\
\Pi_{j=0}^{m-3}(s_{n+1-m}-j) & ... & \Pi_{j=0}^{m-3}(s_{n-1}-j)
\end{array}
\right )
$$
Therefore
$$D=\Pi_{n+1-m\le i<j\le n}(s_i-s_j). $$
The second term comes from $c_{i,i-1}$ and its conjugate where
$$i\in I_2:=\{j\in I;\quad \gamma_j=0. \},\quad I_1=\{i;\quad i\not \in I_2\}. $$
For $i\not \in I_2$, the corresponding term with $c_{i,i-1}$ is not useful since $s_i-s_{i-1}>1$. Also $c_{i,i-1}=0$ if $\gamma_i$ is not an integer. In the expression of $e^{-U^m}$ we see that if $i>n+1-m$ and $i\in I_2$, there is no contribution from
$c_{i,i-1}$ in this case since the second term in each entry of this column (with $c_{i,i-1}$) is just the $c_{i,i-1}$ multiple of the major term the next column (by $s_i=s_{i-1}+1$). Thus the calculation involving this $c_{i,i-1}$ is zero. We only need to consider the following case: $n+1-m\in I_2$ and the coefficient of $c_{n+1-m,n-m}$ is  the determinant of the following matrix:
$$\left(\begin{array}{cccc}
z^{s_{n-m}}& z^{s_{n+2-m}}&...& z^{s_n}\\
s_{n-m}z^{s_{n-m}-1}& s_{n+2-m}z^{s_{n+2-m}-1}&...& s_n z^{s_n-1}\\
\vdots & \vdots & ... & \vdots \\
\Pi_{j=0}^{m-2}(s_{n-m}-j)z^{s_{n-m}+1-m} & \Pi_{j=0}^{m-2}(s_{n+2-m}-j)z^{s_{n+2-m}+1-m}&... & \Pi_{j=0}^{m-2}(s_n-j)s^{s_n+1-m}
\end{array}
\right )
$$
Note that the powers of $z$ in the first column are $2$ less than those in the second column.  Let
$$ D_1=det\left (\begin{array}{cccc}
1 & 1 & ... & 1\\
s_{n-m} & s_{n+2-m} & ... & s_n \\
\vdots & \vdots & ... & \vdots \\
\Pi_{j=0}^{m-2}(s_{n-m}-j) & \Pi_{j=0}^{m-2}(s_{n+2-m}-j) & ... & \Pi_{j=0}^{m-2}(s_n-j)
\end{array}
\right )
$$
It is easy to evaluate $D_1$ in the same way that $D$ was evaluated:
$$D_1=\Pi_{i\ge n+2-m}(s_i-s_{n-m})\Pi_{n+2-m\le i<j\le n}(s_j-s_i). $$
If $n+1-m\in I_1$,
$$det(B_{\lambda})=\lambda_{n+1-m}...\lambda_nDz^{S_m}(1+O(1/z^2)). $$
If $n+1-m\in I_2$,
$$det(B_{\lambda})=\lambda_{n+1-m}...\lambda_nDz^{S_m}(1+\frac{D_1}Dc_{n+1-m,n-m}\frac 1{z}+O(1/z^2)). $$
Correspondingly
$$det(\bar B)=\left\{\begin{array}{ll}
D\bar z^{S_m}(1+O(1/\bar z^2)),\quad \mbox{if}\quad n+1-m\in I_1,\\
\\
D\bar z^{S_m}(1+\frac{D_1}D\bar c_{n+1-m,n-m}\frac 1{\bar z}+O(\frac 1{\bar z^2})),\quad \mbox{if}\quad n+1-m\in I_2.
\end{array}
\right.
$$
Consequently if $n+1-m\in I_1$,
\begin{equation}\label{i1-1}
e^{-U^m}=\lambda_{n+1-m}...\lambda_n|z|^{2S_m}D^2(1+O(1/|z|^2)).
\end{equation}
If $n+1-m\in I_2$,
\begin{align} \label{exp-u}
&e^{-U^m}\\
=&\lambda_{n+1-m}..\lambda_nD^2|z|^{2S_m}(1+\frac{D_1}{D}(\frac{c_{n+1-m}}z+\frac{\bar c_{n+1-m}}{\bar z})+O(1/|z|^2))\nonumber\\
&=\lambda_{n+1-m}..\lambda_nD^2|z|^{2S_m}(1+\frac{2D_1}{D}(\alpha_{n+1-m}\cos\theta+\beta_{n+1-m}\sin\theta)r^{-1}+O(r^{-2})) \nonumber
\end{align}
where $c_{n+1-m}=\alpha_{n+1-m}+\sqrt{-1}\beta_{n+1-m}$, $r=|z|$.
Thus by (\ref{i1-1}) and (\ref{exp-u}) if $n+1-m\in I_2$
\begin{equation}\label{U-diff}
\left\{\begin{array}{ll}
\displaystyle{-\frac{\partial U^l}{\partial \alpha_{n+1-m}}=\frac{2D_1}D\frac{\delta_{lm}\cos \theta}{r}+O(1/r^2),}, \quad \delta_{lm}=\left\{\begin{array}{ll} 1, \, l=m\\
0,\,\, l\neq m.
\end{array}
\right. \\
\\
\displaystyle{-\frac{\partial U^l}{\partial \beta_{n+1-m}}=\frac{2D_1}D\frac{\delta_{lm}\sin \theta}{r}+O(1/r^2).}
\end{array}
\right.
\end{equation}

\section{Proof of Theorem \ref{local-est}}

Recall that $\tilde v_k=(\tilde v_{k,1},..,\tilde v_{k,n})$ is defined in (\ref{vik-d}) and it satisfies (\ref{vik}). The main assumption of this article is that
$\tilde v_k$ converges to a global $SU(n+1)$ Toda system after scaling ( see (\ref{maj-a})).

In the first step of the proof we invoke the main result in \cite{lin-wei-zhang-jems}. There exists a sequence of global solutions $\tilde U_k=(\tilde U_{k,1},..,\tilde U_{k,n})$ of
$$\Delta \tilde U_{k,i}(y)+\sum_{j=1}^n a_{ij}|y|^{2\gamma_j} h_{k,j}(0)e^{\tilde U_{k,j}(y)}=0,\quad \mbox{ in }\quad \mathbb R^2, \quad i=1,...,n$$
such that the following holds:
\begin{enumerate}
\item Let $\lambda_{k,i}$ ($i=0,...,n$), $c_{k,ij}$ ($0\le i<j\le n$) be the parameters in the definition of $\tilde U_k$ and let $\lambda_i$ ($i=0,..,n$) and $c_{ij}$ ($0\le i<j\le n$) be the parameters in the definition of $\tilde v$ in (\ref{maj-a}). Then along a subsequence $\lambda_{k,i}\to \lambda_i$ ($i=0,...,n$) and
    $c_{k,ij}\to c_{ij}$. As a result $\tilde U_{k,i}$ converges to $\tilde v_i$ uniformly over any fixed compact subset of $\mathbb R^2$. Here we use the harmless assumption $\lim_{k\to \infty}h_{k,i}(0)= 1$.
\item There exist distinct points $p_1,...,p_l\in \mathbb R^2$ with $l\le n^2+2n$ such that
$$\tilde v_{k,1}(p_m)-\psi_{k,1}(p_m)=\tilde U_{k,1}(p_m),\quad 1\le m\le l. $$
In other words the first component of $\tilde v_k-\psi_k$ and $\tilde U_k$ agree at $l$ points. In \cite{lin-wei-zhang-jems} these points are determined in a way that $c_{k,ij}$ and $\lambda_{k,ij}$ do not tend to infinity.
\item Let
$$w_{k,i}(y)=\tilde v_{k,i}(y)-\tilde U_{k,i}(y)-\psi_{k,i}(\epsilon_k y)
\quad \mbox{in}\quad \Omega_k:=B(0,\epsilon_k^{-1}). $$
It holds:
\begin{equation}\label{crude-w}
|w_{k,i}(y)|\le C\epsilon_k(1+|y|),\quad y\in \Omega_k.
\end{equation}
\end{enumerate}
It is already established in \cite{lin-wei-ye} that global solutions $U=(U_1,..,U_i)$ of (\ref{toda-nk}) satisfies
$$U_i(y)=(-4-2\gamma_{n+1-i})\log |y|+O(1), \quad |y|>1. $$
What we need is a little more specific expansion of $\tilde U_{k,i}$:
\begin{equation}\label{tUik}
\tilde U_{k,i}(y)=(-4-2\gamma_{n+1-i}-2\gamma_i)\log |y|+c_{k,i}+O(1/|y|),
\end{equation}
for $|y|>1$ and $i=1,..,n$, where $c_{k,i}$ are uniformly bounded.
To see why (\ref{tUik}) holds, we let
$$\hat U_{k,i}(y)=\tilde U_{k,i}(y)+(4+2\gamma_{n+1-i}+2\gamma_i)\log |y|, \quad |y|>1. $$
Then it is easy to see that $\hat U_{k,i}$ satisfies
$$\Delta \hat U_{k,i}(y)+\sum_j a_{ij} |y|^{-2\gamma_{n+1-j}-4}e^{\hat U_{k,j}}=0, \quad |y|>1 $$
and $\hat U_{k,i}$ is bounded at infinity because all the parameters in the definition of $U_k$ are bounded. Making a Kelvin transformation of $\hat U_{k,i}$:
$$V_{k,i}(z)=\hat U_{k,i}(\frac{z}{|z|^2}), \quad |z|<\frac 12,  $$
we have
$$\Delta V_{k,i}(z)+\sum_j |z|^{2\gamma_{n+1-j}}e^{V_{k,j}(z)}=0, \quad |z|<1. $$
Since $V_{k,i}$ is bounded around $0$, $\gamma_i\ge 0$ for all $i$, $V_{k,i}\in C^{1, \alpha}(B_{1/2})$ for all $\alpha\in (0,1)$. From the expansion of $V_{k,i}$ near $0$ we see that (\ref{tUik}) holds.

The equation for $w_{k,i}$ is
\begin{align}\label{wik}
&\Delta w_{k,i}(y)+\sum_j a_{ij}|y|^{2\gamma_j}h_{k,j}(0)e^{\xi_j^k}w_{k,j}\\
=&\sum_j a_{ij}|y|^{2\gamma_j}(h_{k,j}(0)-h_{k,j}(\epsilon_ky)e^{\psi_{k,j}(\epsilon_ky)})e^{\tilde U_{k,j}},\quad i=1,..,n \nonumber
\end{align}
where
$h_{k,i}(0)e^{\xi_i^k(x)}$ is obtained by mean value theorem.
Let
$$w^i_k=\sum_ja^{ij}w_{k,j}$$  then
(\ref{wik}) becomes
\begin{align}\label{wik-1}
&\Delta w^i_k+|y|^{2\gamma_i}h_{k,i}(0)e^{\xi_i^k}w_{k,i}\\
=& |y|^{2\gamma_i}(h_{k,i}(0)-h_{k,i}(\epsilon_ky)e^{\psi_{k,i}(\epsilon_ky)})e^{\tilde U_{k,i}},\quad \mbox{ in }\quad \Omega_k. \nonumber
\end{align}
Let $\phi_k=(\phi_{k,1},...,\phi_{k,n})$ be solutions of the linearized system
\begin{equation}\label{phi-do}
\Delta \phi_{k,i}+\sum_j a_{ij}|y|^{2\gamma_j} h_{k,j}(0)e^{\tilde U_{k,j}}\phi_{k,j}=0 \quad \mbox{in} \quad \mathbb R^2.
\end{equation}
Set
$$\phi^i_k=\sum_ja^{ij}\phi_{k,j}. $$
 Then
 \begin{equation}\label{phi-up}
\Delta \phi^i_k+|y|^{2\gamma_i}h_{k,i}(0)e^{\tilde U_{k,i}}\phi_{k,i}=0,\quad \mbox{ in }\quad \mathbb R^2.
\end{equation}

Let  $\tilde \Omega_k=B(0,\frac 12\epsilon_k^{-1})$.
Multiplying $\phi_{k,i}$ to both sides of (\ref{wik-1}) and taking the summation on $i$, we have,
\begin{align}\label{wik-11}
&\int_{\tilde \Omega_k}\sum_i\Delta (\sum_{j}a^{ij}w_{k,j})\phi_{k,i}+\int_{\tilde \Omega_k}\sum_i|y|^{2\gamma_i}h_{k,i}(0)e^{\xi_i^k}w_{k,i}\phi_{k,i}\\
=& \sum_i\int_{\tilde \Omega_k} |y|^{2\gamma_i}\bigg (h_{k,i}(0)-h_{k,i}(\epsilon_ky)e^{\psi_{k,i}(\epsilon_ky)}\bigg )e^{\tilde U_{k,i}}\phi_{k,i}, \nonumber
\end{align}

Using integration by parts, we write the first term on the left hand side of (\ref{wik-11}) as
\begin{align*}
&\int_{\tilde \Omega_k}\Delta(\sum_{i,j}a^{ij}w_{k,j})\phi_{k,i}\\
=&\int_{\partial \tilde \Omega_k}(\sum_{i,j} a^{ij}\partial_{\nu}w_{k,j}\phi_{k,i}-\sum_{i,j} a^{ij}w_{k,j}\partial_{\nu}\phi_{k,i})dS+\int_{\tilde \Omega_k}\sum_{i,j} a^{ij}w_{k,j}\Delta \phi_{k,i}. \\
=&\int_{\partial \tilde \Omega_k}\sum_i (\partial_{\nu}w^i_k\phi_{k,i}-w^i_k\partial_{\nu}\phi_{k,i})dS+\int_{\tilde \Omega_k} \sum_i w_{k,i}\Delta \phi^i_k.
\end{align*}

The left hand side of (\ref{wik-11}) now becomes
\begin{equation}\label{wik-lhs}
\int_{\partial \tilde \Omega_k}\sum_i(\partial_{\nu}w^i_k\phi_{k,i}-w^i_k\partial_{\nu}\phi_{k,i})dS+\int_{\tilde \Omega_k}\sum_i(\Delta \phi^i_k+|y|^{2\gamma_i}h_{k,i}(0)e^{\xi_i^k}\phi_{k,i})w_{k,i}dy.
\end{equation}
By (\ref{crude-w}) and (\ref{tUik})
$$h_{k,i}(0)e^{\xi_i^k}-h_{k,i}(0)e^{\tilde U_{k,i}}=O(\epsilon_k)(1+|y|)^{-3-2\gamma_i-2\gamma_{n+1-i}} $$
$\phi_{k,i}$ will be chosen to satisfy
\begin{equation}\label{phi-ik-cr}
\left\{\begin{array}{ll} \phi_{k,i}(x)=(d_{k,i}\cos \theta +q_{k,i}\sin \theta)/r+O(1/r^2)\\
\\
\partial_r\phi_{k,i}(x)=-\displaystyle{\frac{d_{k,i}\cos\theta+q_{k,i}\sin\theta}{r^2}}+O(1/r^3),
\end{array}
\right.
\quad r=|x|>1,
\end{equation}
where $d_{k,i}$ and $q_{k,i}$ are bounded sequences of constants.
Thus by (\ref{crude-w}), (\ref{tUik}) and the estimate of $\phi_{k,i}$  in (\ref{phi-ik-cr}) above, we have
\begin{equation}\label{2-term}
|y|^{2\gamma_i}e^{\xi_i^k}\phi_{k,i}w_{k,i}-|y|^{2\gamma_i}e^{\tilde U_{k,i}}\phi_{k,i}w_{k,i}=O(\epsilon_k^2)(1+|y|)^{-3-2\gamma_{n+1-i}},\quad |y|>1.
\end{equation}
Using (\ref{2-term}) and (\ref{phi-up}) in the second term of (\ref{wik-lhs}) we obtain the following estimate easily
$$\int_{\tilde \Omega_k}\sum_i(\Delta \phi^i_k+|y|^{2\gamma_i}h_{k,i}(0)e^{\xi_j^k}\phi_{k,i})w_{k,i}dy=O(\epsilon_k^2). $$
We further claim that the first term of (\ref{wik-lhs}) is $O(\epsilon_k^2)$, which follows immediate from (\ref{phi-ik-cr}) and the following estimate:
\begin{equation}\label{wik-2}
\nabla w_{k,i}=O(\epsilon_k^2),\quad w_{k,i}=c_{k,i}+O(\epsilon_k),\quad \mbox{ on }\quad \partial \tilde \Omega_k.
\end{equation}
In order to prove (\ref{wik-2}) we first observe that
$$w_{k,i}=c_{k,i}+O(\epsilon_k) \quad \mbox{ on }\quad \partial \Omega_k $$
since $\tilde v_{k,i}-\psi_{k,i}$ is constant on $\partial \Omega_k$ and $\tilde U_{k,i}$ satisfies (\ref{tUik}).
Now we use the Green's representation of $w_{k,i}$ on $\Omega_k$:
\begin{equation}\label{wik-e}
w_{k,i}(y)=\int_{\Omega_k}G_k(y,\eta)(-\Delta w_{k,i}(\eta))d\eta-\int_{\partial \Omega_k}\partial_{\nu}G(y,\eta)w_{k,i}(\eta)dS_{\eta},
\end{equation}
where
\begin{equation}\label{g-exp}
G_k(y,\eta)=-\frac 1{2\pi}\log |y-\eta |+\frac 1{2\pi}\log (\frac{|y|}{\epsilon_k^{-1}}|\frac{\epsilon_k^{-2}y}{|y|^2}-\eta |),\quad y,\eta \in \Omega_k.
\end{equation}
It is easy to see that the second term of (\ref{wik-e}) is a harmonic function with $O(\epsilon_k)$ perturbation on $\partial \Omega_k$. Thus the gradient of this term on $\partial \tilde \Omega_k$ is $O(\epsilon_k^2)$ because the distance from $\partial \tilde \Omega_k$ to $\partial \Omega_k$ is comparable to $\epsilon_k^{-1}$. Therefore in order to prove
(\ref{wik-2}) it suffices to show
\begin{align}\label{wik-3}
&\int_{\Omega_k}\nabla_y G_k(y,\eta)\bigg (\sum_j a_{ij}|\eta |^{2\gamma_j}h_{k,j}(0)e^{\xi_j^k}w_{k,j}(\eta) \\
&\quad -\sum_j a_{ij}|\eta |^{2\gamma_j}(h_{k,j}(0)-h_{k,j}(\epsilon_k \eta)e^{\psi_{k,j}(\epsilon_k \eta)})e^{U_{k,j}}\bigg ) \nonumber\\
=&O(\epsilon_k^2),\quad y\in \Omega_k, \,\, |y|>\frac 14\epsilon_k^{-1}. \nonumber
\end{align}

The proof of (\ref{wik-3}) follows from (\ref{g-exp}), (\ref{crude-w}) and (\ref{tUik}) by standard estimate, so we omit the details. Hence we have established the first estimate of (\ref{wik-2}). To prove the second estimate of (\ref{wik-2}) we just need to show that the oscillation of $w_{k,i}$ on $\Omega_k\setminus \tilde \Omega_k$ is $O(\epsilon_k)$. Let $y_1,y_2\in \Omega_k\setminus \tilde \Omega_k$, using (\ref{wik-e}) we have
\begin{align}\label{osc-2}
&w_{k,i}(y_1)-w_{k,i}(y_2)\\
=&\int_{\Omega_k}(G_k(y_1,\eta)-G(y_2,\eta))\bigg (\sum_j a_{ij}|\eta |^{2\gamma_j}h_{k,j}(0)e^{\xi_j^k}w_{k,j}(\eta)\nonumber\\
&\quad -\sum_j a_{ij}|\eta |^{2\gamma_j}(h_{k,j}(0)-h_{k,j}(\epsilon_k \eta)e^{\psi_i^k(\epsilon_k \eta)})e^{U_{k,j}}\bigg )\nonumber\\
&-\int_{\partial \Omega_k}(\partial_{\nu}G_k(y_1,\eta)-\partial_{\nu}G_k(y_2,\eta))w_{k,i}(\eta)dS_{\eta}. \nonumber
\end{align}
The last term in (\ref{osc-2}) is $O(\epsilon_k)$ because it is the difference of of two points of a harmonic function whose oscillation on $\partial \Omega_k$ is $O(\epsilon_k)$. Writing $G_k(y_1,\eta)-G_k(y_2,\eta)$ as
$$G_k(y_1,\eta)-G_k(y_2,\eta)=\nabla_1G_k(y^*,\eta)\cdot (y_1-y_2) $$
where $\nabla_1$ means differentiation with respect to the first component, $y^*$ is between $y_1$ and $y_2$ ( $y^*\in \Omega_k\setminus \tilde \Omega_k$),
we see that the first term of (\ref{osc-2}) is $O(\epsilon_k)$ by (\ref{wik-3}). Thus (\ref{wik-2}) is established.
By (\ref{wik-2}) the left hand side of (\ref{wik-11}) is $O(\epsilon_k^2)$.
Equation (\ref{wik-11}) now becomes
\begin{equation}\label{wik-rhs}
\sum_i\int_{\tilde \Omega_k}|y|^{2\gamma_i}(h_{k,i}(0)-h_{k,i}(\epsilon_ky)e^{\psi_{k,i}(\epsilon_ky)})e^{\tilde U_{k,i}(y)}\phi_{k,i}(y)dy=O(\epsilon_k^2).
\end{equation}
Let
\begin{equation}\label{bar-hik}
\bar h_{i,k}=h_{k,i}e^{\psi_{k,i}},
\end{equation}
by  (\ref{tUik}) and (\ref{phi-ik-cr}) we have
\begin{align*}
&|y|^{2\gamma_i}(h_{k,i}(0)-h_{k,i}(\epsilon_ky)e^{\psi_{k,i}(\epsilon_ky)})e^{\tilde U_{k,i}(y)}\phi_{k,i}(y)\\
=&\epsilon_k (\partial_1\bar h_{k,i}(0)y_1+\partial_2 \bar h_{k,i}(0)y_2)|y|^{2\gamma_i}h_{k,i}(0)e^{\tilde U_{k,i}}\phi_{k,i}+O(\epsilon_k^2)(1+|y|)^{-3-2\gamma_{n+1-i}}.
\end{align*}
Therefore (\ref{wik-rhs}) is reduced to
$$\sum_i\int_{\tilde \Omega_k}(\partial_1\bar h_{k,i}(0)y_1+\partial_2 \bar h_{k,i}(0)y_2)|y|^{2\gamma_i}h_{k,i}(0)e^{\tilde U_{k,i}}\phi_{k,i}=O(\epsilon_k). $$
Using (\ref{phi-do}) and integration by parts we have
\begin{align}\label{eq-ob}
&\sum_i\int_{\tilde \Omega_k}(\partial_1\bar h_{k,i}(0)y_1+\partial_2 \bar h_{k,i}(0)y_2)|y|^{2\gamma_i}h_{k,i}(0)e^{\tilde U_{k,i}}\phi_{k,i}\\
=&\sum_i\int_{\partial \tilde \Omega_k}(y_1\partial_1\bar h_{k,i}(0)+y_2\partial_2\bar h_{k,i}(0))(\frac{\phi^i_k}{|y|}-\partial_{\nu}\phi^i_k)\nonumber\\
=&\sum_i\int_0^{2\pi}(\cos \theta \partial_1 \bar h_{k,i}(0)+\sin\theta \partial_2 \bar h_{k,i}(0))(\phi^i_k-r\partial_r\phi^i_k)rd\theta \nonumber\\
=&O(\epsilon_k)\nonumber
\end{align}
From (\ref{bar-hik}) we have
\begin{equation}\label{g-bh}
\nabla \bar h_{k,i}(0)=(\nabla \log h_{k,i}(0)+\nabla \psi_{k,i}(0))h_{k,i}(0).
\end{equation}
If $$\phi^i_k=\frac 1r(d_{k,i}\cos\theta+q_{k,i}\sin \theta)+O(1/r^2) $$
we obtain from (\ref{eq-ob})
\begin{align}\label{van-1}
&\sum_i\int_0^{2\pi}(\cos\theta \partial_1 \bar h_{k,i}(0)+\sin\theta \partial_2 \bar h_{k,i}(0))r(\phi^i_k-r\partial_r\phi^i_k)d\theta\\
=&\pi \sum_i(d_{k,i}\partial_1\bar h_{k,i}(0)+q_{k,i}\partial_2 \bar h_{k,i}(0))+O(\epsilon_k)\nonumber
\end{align}
where $r=\epsilon_k^{-1}/2$.

Recall that $(U^1_k,...,U^n_k)$ are described by up to $n^2+2n$ parameters. In particular we write
$$c_{k,n+1-i}=\alpha_{k,n+1-i}+\sqrt{-1}\beta_{k,n+1-i}. $$
 For each $i$ such that $n+1-i\in I_2$, we differentiate the real and imaginary parts of $(U_k^1,...,U_k^n)$ to get (by (\ref{U-diff}))
$$\Phi_{\alpha_{n+1-i}}
=\left(\begin{array}{c}
-\frac{\partial U_k^{1}}{\partial \alpha_{k,n+1-i}}\\
\vdots\\
-\frac{\partial U_k^{n}}{\partial \alpha_{k,n+1-i}}
\end{array}
\right )=
\left(\begin{array}{c}
\frac{2D_1}D\delta_{1,i}\cos\theta/r\\
\vdots\\
\frac{2D_1}D\delta_{n,i}\cos\theta/r
\end{array}
\right)+O(\frac 1{r^2})
$$
and
$$\Phi_{\beta_{n+1-i}}=\left(\begin{array}{c}
-\frac{\partial U_k^{1}}{\partial \beta_{k,n+1-i}}\\
\vdots\\
-\frac{\partial U_k^{n}}{\partial \beta_{k,n+1-i}}
\end{array}
\right )=
\left(\begin{array}{c}
\frac{2D_1}D\delta_{1,i}\sin\theta/r\\
\vdots\\
\frac{2D_1}D\delta_{n,i}\sin\theta/r
\end{array}
\right)+O(\frac 1{r^2})
$$
For $n+1-i\in I_2$, letting $(\phi_k^1,...,\phi_k^n)'$ be $\Phi_{\alpha_{n+1-i}}$ and $\Phi_{\beta_{n+1-i}}$ respectively in (\ref{van-1}) we have
$$|\nabla \bar h_{k,i}(0)|=O(\epsilon_k),\quad \mbox{ if}\quad n+1-i\in I_2 $$
which is, by (\ref{g-bh}),
\begin{equation}\label{van-2}
\nabla (\log h_{k,i})(0)+\nabla \psi_{k,i}(0)=O(\epsilon_k),\quad \mbox{if}\quad n+1-i\in I_2.
\end{equation}

Theorem \ref{local-est} is established. $\Box$

\begin{rem} In \cite{lin-wei-zhang-adv} a $\partial_z^2$ estimate was established for fully bubbling solutions of nonsingular Toda systems. In this article
it does not seem to be possible to get the $\partial_z^2$ estimates for the nonsingular $SU(n+1)$ Toda system. The reason is if the estimate of $\nabla h_i^k(0)$ is not obtained for $n+1-i\in I_1$, the corresponding $w_{k,i}$ cannot be improved to $O(\epsilon_k^2)$ over compact subsets of $\mathbb R^2$. This lack of accuracy prevents us from getting estimates on the second derivatives of $h_i^k$.
\end{rem}

\end{document}